\documentstyle[12pt]{amsart}
\theoremstyle{plain}
\newtheorem{Thm}{}[section]

\title {Existence of coherent systems.}
\author{ Montserrat Teixidor i Bigas}
\address{ Mathematics Department, Tufts University, Medford MA
02180, U.S.A.}

\begin{document}
\maketitle
\begin{section}{Introduction}

Let $C$ be a projective, algebraic non-singular curve of given
genus $g$. A coherent system on $C$ is a pair  $(E,V)$ where $E$
is a vector bundle and $V$ is a subspace of the space of sections
of $E$.
  A coherent subsystem is a pair $(E',V')$ where $E'$ is a
subbundle of $E$ and $V'\subset V\cap H^0(E')$.

Given a real number $\alpha >0$ a coherent system is said to be
$\alpha$-(semi)-stable if for every proper coherent subsystem
$${{\deg }\ E'+\alpha \ {\dim}\ V'\over {\mbox {\rom rk} }\ E'}<(\leq )
{{\deg }\ E+\alpha \ {\dim}\ V\over {\mbox {\rom rk }}E}.$$

Fix values of the degree $d$, the dimension $k$ of $V$ and the
rank $r$ of $E$ as well as a rational number $\alpha >0$. The set
of coherent systems that are $\alpha $-stable can be given the
structure of a moduli space (see \cite{KN} \cite{LP}). This moduli
space, if non-empty, has dimension at each point at least equal to
the Brill-Noether number (see \cite{BGMN})
$$\rho =r^2(g-1)+1-k(k-d+r(g-1)).$$
In fact, it is expected that under mild additional conditions, the
moduli space will be non-empty of dimension $\rho$ provided that
$\rho$ is positive.

Coherent systems have received a great deal of attention in the
last few years. For an (almost up to date) detailed exposition of
the current knowledge, the reader is advised to read the
introduction in \cite{BGMMN}.

 For genus one, these spaces have been described precisely in \cite{LN}. For higher genus,
 the most pressing question is to have conditions that ensure the non-emptiness of these moduli spaces.
  Such conditions have been given for $k\le r$
  in \cite{BG, BGMN, BGMMN}. A great deal is known for $k=r+1$
  (see \cite{BGMN,B}). In this paper we shall concentrate on the
  case $k>r$ for generic curve

\begin{Thm} {\bf Theorem}
Let $C$ be a generic non-singular curve of genus $g\ge 2$.
 Let $d,r,k$ be positive integers with  $ k>r$. Write
$$d=rd_1+d_2,\ k=rk_1+k_2, \ d_2<r,\ k_2<r$$ and all $d_i,\ k_i$
non-negative integers. Assume that
$$(*) g-(k_1+1)(g-d_1+k_1-1)\ge 1 ,\ 0\not= d_2\ge k_2 $$
$$(**) g-k_1(g-d_1+k_1-1)> 1 ,\ d_2= k_2=0$$
$$(***) g-(k_1+1)(g-d_1+k_1)\ge 1 ,\ d_2< k_2.$$
Then, for any positive value of $\alpha$ the moduli space of
coherent systems of rank $r$ degree $d$ and with $k$ sections on
$C$ is non-empty and has one component of the expected dimension
$\rho$.
\end{Thm}

\begin{Thm} {\bf Remark}
The bounds in the Theorem are those in \cite{T2} and are not meant
to be the best possible. The methods in this paper can be used to
improve them although it is likely that this will require a case
by case analysis (see for instance \cite{T1} for a similar problem
in the case of rank two).
\end{Thm}

Acknowledgements: The author was partially supported by a
Radcliffe Fellowship during the preparation of this work. She is a
member of the research group VBAC and wants to thank P.Newstead
for organizing a workshop on coherent systems in August 2005. I
would also like to thank the referee for a very careful reading of
the manuscript.
\end{section}

 \begin{section}{ Background on reducible curves}

This moduli space, if non-empty, has dimension at each point at
least equal to the Brill-Noether number (see \cite{BGMN})
$$\rho =r^2(g-1)+1-k(k-d+r(g-1).$$
We shall use a special type of reducible curve that we shall call
a chain of elliptic curves and is defined as follows:

\begin{Thm} {\bf Definition} \label{cce}
Let $C_1...C_g$ be elliptic curves. Let $P_i,Q_i$ be generic
points in $C_i$. Then $C_0$ is the  chain obtained by gluing the
elliptic curves by identifying the point $Q_i$ in $C_i$
 to the point $P_{i+1}$ in $C_{i+1},\ i=1...g-1$. When there is no
 danger of confusion, we shall write $P,Q$ instead of $P_i,Q_i$.
\end{Thm}

 We shall be using the following well-known fact

\begin{Thm} \label{remark}
{\bf Lemma} Let $C$ be an elliptic curve and $L$ a line bundle of
degree $d$ on $C$. One can define a subspace of dimension $k$ of
sections of  $L$ by specifying the $k$ distinct (minimum) desired
vanishings of a basis of the subspace at two different points
$P,Q$ so that the sum of the
 corresponding vanishings at $P$ and $Q$ is $d-1$.
In the case when $L={\cal O}(aP+(d-a)Q)$ two of the vanishings
could be chosen to be $a,\ d-a$ adding to $d$ rather than $d-1$.
These are the
 only two vanishings that can add up to $d$ if $P,Q$ are generic.

Similarly, let $E$ be a vector bundle obtained as the sum of $r$
line bundles of degree $d$. Then one can find a space $V$ of
dimension $k$ of sections of $E$ with desired vanishings
$a_1,...,a_k$ and $b_1,...,b_k$ respectively at two points $P,Q$
if  $a_i+b_{k-i}\le d-1$ and each possible value of $a_i,b_i$
appears at most $r$ times.

If $E$ is an indecomposable vector bundle of rank $r$ and degree
$rd_1+d_2$ with $0<d_2<r$, there are no sections whose orders of
vanishing at two points add up to a number greater than $d_1$. The
space of sections vanishing to order $a$ at $P$ and $d_1-a$ at $Q$
has dimension $d_2$.

\end{Thm}
The proof is left to the reader (or see \cite{T1})

When dealing with reducible curves, the notion of a line bundle
and a space of its sections needs to be replaced by the analogous
concept of limit linear series as introduced by Eisenbud and
Harris.
 A similar definition can be given
for vector bundles (cf\cite{T2} ). For the convenience of the
reader, we reproduce this definition here

\begin{Thm}\label{lls}
 {\bf Limit linear series}
A limit linear series of rank $r$, degree $d$ and dimension $k$ on
a chain of $M$ (not necessarily elliptic) curves consists of data
I,II below for which data III, IV exist satisfying conditions
a)-c)

I) For every component $C_i$, a vector bundle $E_i$ of rank $r$
and degree $D_i$ and a $k$-dimensional space $V_i$ of sections of
$E_i$

II) For every node obtained by gluing $Q_i$ and $P_{i+1}$, an
isomorphism of the projectivisation of the fibers $(E_i)_{Q_i}$
and $(E_{i+1})_{P_{i+1}}$

III) A positive integer $b$

IV) For every node obtained by gluing $Q_i$ and $P_{i+1}$ bases
$s^t_{Q_i}, s^t_{P_{i+1}},\ t=1...k$ of the vector spaces $V_i,
V_{i+1}$ of I.

\medskip
Subject to the conditions

a) $\sum_{i=1}^M D_i-r(M-1)b=d$

b) The orders of vanishing at $Q_i, P_{i+1}$ of the sections of
the chosen bases satisfy
$ord_{Q_i}s^t_{Q_i}+ord_{P_{i+1}}s^t_{P_{i+1}}\ge b$

c) Sections of the vector bundles $E_i(-bP_i), E_i(-bQ_i)$ are
completely determined by their values at the nodes.
\end{Thm}
\end{section}

\begin{section}{Sketch of the proof}

We now give a general outline of the proof of the theorem.

 We want to show that for any positive $\alpha$ and given degree $d$,
 rank $r$ and number of sections $k$, the set of $\alpha$-stable
coherent systems satisfying the conditions of the Theorem is
non-empty. We shall choose $E$ to be a stable vector bundle (with
the usual meaning for stability) and $V$ a generic subspace of
dimension $k$ of its sections. Because of the stability of $E$,
for every subbundle $E'$ of $E$ one has the inequality
$${{\deg} E'\over {\mbox{\rom rk}}E'}<{{\deg} E\over {\mbox{\rom rk}}E}.$$
Therefore, in order to prove the stability of the coherent system,
it suffices to show that for every subbundle $E'$ of $E$, the
dimension $k'$ of the vector space $H^0(E')\cap V$ satisfies
$${k'\over {\mbox {\rom rk}}E'}\leq {k\over {\mbox{\rom rk}}E}.$$
In fact, if $E$ is stable then $(E,V)$ is $\alpha$-stable for
arbitrarily large $\alpha$ if and only if the condition above is
satisfied.

We need to prove the result for the generic curve. It suffices to
prove it for a particular curve and a coherent system which is a
member of a family of dimension precisely $\rho$ on that curve.
This follows from the fact that the dimension of a family of
coherent systems on a family of curves ${\mathcal C}\rightarrow B$
is at least $\dim B+\rho $. If the dimension of one of the fibers
over $B$ is precisely $\rho$, then the family projects onto the
base (see for instance \cite{ S, T2, T1} for similar arguments).

The curve we shall use will be a singular curve of the special
type that we described \ref{cce}.
 We then construct a limit coherent system on the chain of elliptic curves  and
show  as outlined above that no subbundle has enough sections to
contradict stability.
\end{section}

\begin{section}{The limit linear series}

  As in the statement of the Theorem, write $d=rd_1+d_2, k=rk_1+k_2$. Let $h$ be the greatest
common divisor of $d$ and $r$ and $d=\bar dh,\ r=\bar rh$. We
describe next limit linear series on $C_0$, we need to distinguish
between various cases depending on the relative values of
$d_2,k_2$.  The number $b$ that appear in III of \ref{lls}, will
be taken to be $d_1$ in all cases.

Assume first $d_2< k_2$ and $g-d_1+k_1\ge 1$. The latter condition
is equivalent to the fact that a generic vector bundle of rank $r$
and degree $d$ has fewer than $k$ sections.

We start by describing the vector bundles and gluings that appear
in the definition of limit linear series. We leave the description
of the spaces of sections to later.

 On $C_1$ take the vector bundle to be the direct sum of $h$
 generic vector bundles of rank $\bar r$ and degree $\bar d$.

 On the curve $C_i,\ i=2,...,k_1+2$ take
the vector bundle
$$ {\mathcal O}((2i-3)P+(d_1-2i+3)Q)^{k_2-d_2}\oplus
L^i_1\oplus...\oplus L^i_{r-k_2+d_2}$$ where the $L^i_j$ are
generic line bundles of degree $d_1$. Consider the space $\bar
W_1$ of dimension $d_2$ of sections of $E_1$ with maximum
vanishing $d_1$ at $Q_1$. Take the gluing of $E_{1,Q_1}$ with
$E_{2,P_2}$ so that the image of  $W_{1Q_1}$ by the map
$H^0(E_1(-d_1Q_1)\otimes {\mathcal O}_C\rightarrow
E_{1Q_1}(-d_1Q_1)$ lies
 inside the space of dimension $r-k_2+d_2$ of the fiber
of $E_2$ at $P_2$ generated by the subbundle $L^2_1\oplus...\oplus
L^2_{r-k_2+d_2}$.

Take the gluing at the next $k_1$ nodes so that the subbundle
$L^i_1\oplus...\oplus L^i_{r-k_2+d_2}$ on the curve $C_i$ glues
with the analogous subbundle on the curve $C_{i+1}$.

Consider inside the space of sections of  $E_2(-(d_1-1)Q_2)$ the
subspace $\bar W_2$ of those gluing at $P_2$ with $\bar W_1$. Let
$ \tilde W_2$ be the subspace of sections of $
E_2(-(d_1-1)Q_2-P_2)$. Define $W_2= \bar W_2+\tilde W_2$.

 For $i\ge 3$, define $W_i$ by induction in the following way:
Consider the space of sections of  $E_i(-(i-2)P_i-(d_1-i+1)Q_i)$.
It has dimension $r$ and maps bijectively on the fibers of $P_i$
and $Q_i$. Take as $W_i$ the subspace of dimension $k_2$ that
glues with $W_{i-1}$ at $P_i$.

Note that by our assumption (***) $(k_1+1)(g+k_1-d_1)+1\le g$.  On
the curves
$$C_{(k_1+1)\alpha +\beta +1 },\ \alpha =1,...,g+k_1-d_1-1, \
\beta =1,...,k_1+1$$ take the following vector bundles:

 If $\beta =1$ take
$$ {\mathcal O}((k_1\alpha +1)P+(d_1-k_1\alpha -1)Q)^{k_2}\oplus
L^i_1\oplus...\oplus L^i_{r-k_2}$$

For $\beta =2,...k_1+1$ take
$$ {\mathcal O}((k_1\alpha +2\beta -1)P+(d_1-k_1\alpha -2\beta
+1)Q)^{r}.$$

For $\beta =2,..., k_1+1$, the gluing is generic.

 For $\beta =1$, $\alpha =1$
take the gluing so that the subbundle
$$ {\mathcal O}((k_1+1)P+(d_1-k_1 -1)Q)^{k_2}$$
 glues with $W_{k_1+2}$.

Define $W_i$ as before where for $i=(k_1+1)\alpha +1+\beta$ the
spaces of sections to be considered are
$$H^0(E_i(-(\alpha k_1+\beta
-1)P -(d_1-\alpha k_1-\beta )Q), \beta \not= 1$$
$$H^0(E_i(-(\alpha k_1+1 )P -(d_1-\alpha k_1-1)Q), \beta
=1.$$

For $\beta =1,\ \alpha >1$ take the gluing so that the fiber at
$Q_{\alpha (k_1+1)+1}$ of the subbundle
$$ {\mathcal O}((k_1\alpha +1)P+(d_1-k_1\alpha -1)Q)^{k_2}$$ glues with $W_{\alpha (k_1+1)+2}$.

On the remaining curves, the vector bundles are direct sums of
generic line bundles of degree $d_1$ and the gluing is generic.

One can still define the subspace $W_i$ as before. Write
$t=g+k_1-d_1$ and use as space of sections
$H^0(E_i(-(i-t-1)P-(d_1-i+t)Q)$.

The resulting vector bundle is stable  because the restriction of
the vector bundle to each component is semistable and the
destabilizing subbundles in the various components do not glue
with each other (see \cite{T3, T4})

We want to compute the dimension of the set of such vector
bundles. We need to add the dimensions of the families of the
restrictions of the vector bundles to each component and the
dimensions of the families of gluings. We then subtract the
dimensions of the families of automorphisms of the vector bundles
restricted to each component and finally add one as the resulting
bundle being stable, it has a one dimensional family of
automorphisms.

The first bundle moves in an $h$ dimensional family and has an $h$
dimensional family of automorphisms. The next $k_1+1$  bundles
move in a $[r-k_2+d_2]$-dimensional family and have a
$[(k_2-d_2)^2+r-k_2+d_2]$-dimensional family of automorphisms. The
first gluing varies in a $[d_2(r-(k_2-d_2))+r(r-d_2)]$-dimensional
family. The next $k_1$ gluings vary in a
$[(r-k_2+d_2)^2+r(k_2-d_2)]$-dimensional family.

 On the curves corresponding to $\beta =1$, the vector bundle
 varies in a $[r-k_2]$-dimensional family and has a
$[(k_2)^2+r-k_2]$-dimensional family of automorphisms while the
gluing depends on $k_2^2+r(r-k_2)$ parameters.

On the curves corresponding to $\beta >1$, the vector bundles are
completely determined, they have an $r^2$-dimensional family of
endomorphisms and the gluings are free.

On the remaining curves, the vector bundles depend on $r$
parameters, have an $r$-dimensional family of endomorphisms and
the gluings are free.

Therefore, the dimension of the family is
$$[h-h]+[r-k_2+d_2-(k_2-d_2)^2-(r-k_2+d_2)
+d_2(r-(k_2-d_2))+r(r-d_2)]$$
$$+k_1[r-k_2+d_2-(k_2-d_2)^2-(r-k_2+d_2)
+(r-(k_2-d_2))^2+r(k_2-d_2)]$$
$$(g+k_1-d_1-1)[r-k_2-k_2^2-(r-k_2)+k_2^2+r(r-k_2)]+(g+k_1-d_1-1)k_1(-r^2+r^2)$$
$$+(g-(k_1+1)(g+k_1-d_1)-1)(r-r+r^2)+1=\rho$$

We need to define the spaces of sections on every component of the
curve and check that this gives rise  to only one limit series for
each bundle as defined above. As pointed out before, we take
$b=d_1$ in III of \ref{lls}.

Similarly to the way we defined $W_i$, we shall define spaces of
sections $W'_i$ this time by descending induction for
$i=k_1+2,...1$. Define first
$$W'_{k_1+2}=H^0(E_{k_1+2}(-(2k_1+1)P-(d_1-2k_1-1)Q)).$$
By definition, $$E_{k_1+2}=({\mathcal
O}((2k_1+1)P+(d_1-2k_1-1)Q))^{k_2-d_2}\oplus L^{k_1+2}_1\oplus
...\oplus L^{k_1+2}_{r-k_2+d_2}$$ where the last $r-k_2+d_2$ terms
are generic line bundles of degree $d_1$. Therefore,
$E_{k_1+2}(-(2k_1+1)P-(d_1-2k_1-1)Q)$ consists of a direct sum of
$k_2-d_2$ copies of the trivial line bundle and $r-k_2+d_2$
generic line bundles of degree zero. It follows that $W'_{k_1+2}$
has dimension $k_2-d_2$.

For $i<k_1+2$, define $W'_i$ by descending induction as the space
of sections of $E_i(-(i+k_1-1)P-(d_1-i-k_1)Q)$ that glues with
$W'_{i+1}$ at $Q_i$.

For $i=1$, the space of sections of the limit linear series is
$$H^0(E_1(-(d_1-k_1)Q_1))\oplus W'_1.$$
The vanishings of the sections of this space at $P_1$ are
$$0,1,k_1-1, k_1$$ where each vanishing is repeated $r$ times
except for the last one which is repeated $k_2$ times. The
vanishings of the sections of this space at $Q_1$ are
$$d_1,d_1-1...,d_1-k_1, d_1-k_1-1$$ where each vanishing is repeated $r$ times
except for the first one which is repeated $d_2$ times and the
last one which is repeated $k_2-d_2$ times.

For $i=2,...,k_1+2$, the space of sections is
$$W_i\oplus H^0(E_i(-(i-1)P-(d_1-2i+3)Q))\oplus
 H^0(L^i_1\oplus...\oplus
L^i_{r+d_2-k_2}(-(2i-3)P-(d_1-2i+2)Q))$$ $$
  \oplus H^0(E_i(-(2i-2)P-(d_1-k_1-i+1)Q))\oplus W'_i.$$

For $i=2$, the vanishings of the sections of this space at $P_2$
are $0,1,2,...,k_1,k_1+1$ where each number is repeated $r$ times
except for the first one which appears $d_2$ times and the last
one that appears $k_2-d_2$ times.

 The vanishings of
the sections of this space at $P_i$ for $i> 2$ are
$$i-2,i-1,.., 2i-4, 2i-3 ...k_1+i-2, k_1+i-1$$ where each vanishing is repeated $r$ times
except for the first one which is repeated $k_2$ times, the
vanishing $2i-4$ which is repeated $r+d_2-k_2$ times and the last
one which is repeated $k_2-d_2$ times.

The vanishings of the sections of this space at $Q_i$ are
$$d_1-i+1,d_1-i...d_1-2i+4, d_1-2i+3, d_1-2i+2...,d_1-k_1-i+1, d_1-k_1-i$$
where each vanishing is repeated $r$ times except for the first
one which is repeated $k_2$ times, the vanishing $d_1-2i+2$ which
is repeated $r+d_2-k_2$ times and the last one which is repeated
$k_2-d_2$ times.

For $i=\alpha (k_1+1)+\beta +1$ the space of sections is
$$W_i \oplus H^0(E_i(-(k_1\alpha +\beta )P-(d_1-k_1\alpha -2\beta +1)Q))\oplus
$$
$$ \oplus H^0(E_i(-(k_1\alpha +2\beta )P-(d_1-k_1(\alpha +1)-\beta -1)Q)).$$
The vanishings of the sections of this space at $P_i$ are
$$i-\alpha-2,i-\alpha-1,..k_1\alpha+2\beta -3,
k_1\alpha+2\beta -1,k_1\alpha+2\beta ,...,k_1(\alpha+1)+\beta
$$ where each vanishing is repeated $r$ times except for the first
one which is repeated $k_2$ times. The vanishings of the sections
of this space at $Q_i$ are
$$d_1-i+\alpha+1,d_1-i+\alpha,..,
d_1-k_1\alpha-2\beta +1,d_1-k_1\alpha-2\beta -1
,...d_1-k_1(\alpha+1)-\beta -1$$ where each vanishing is repeated
$r$ times except for the first one which is repeated $k_2$ times.

Write $t=g+k_1-d_1$. then for $i>t(k_1+1)+1$, the space of
sections is $$W_i\oplus H^0(E_i(-(i-t)P-(d_1-i+t-k_1)Q))$$ The
vanishings of the sections of this space at $P_i$ are
$$i-t-1,i-t,...i-t+k_1-1$$
where each vanishing is repeated $r$ times except for the first
one which is repeated $k_2$ times. The vanishings of the sections
of this space at $Q_i$ are
$$d_1-i+t,d_1-i+t-1,..d_1-i+t-k_1$$
where each vanishing is repeated $r$ times except for the first
one which is repeated $k_2$ times.

One checks then that the vanishings at $P_1$ and $Q_g$ are the
smallest possible, namely $$0,1,...k_1-1,k_1$$ each with
multiplicity $r$ except for $k_1$ which has multiplicity $k_2$.
Moreover, the sum of the vanishings at $Q_i$ and $P_{i+1}$ of
corresponding sections is the minimum possible namely $b=d_1$.
Therefore, deforming one of the $E_i$ at one of the components
$C_i$ to a more general bundle would decrease the vanishing of the
limit linear series at $Q_i$ and would therefore make it
impossible to prolong the limit linear series till $C_g$.
Similarly, deforming one of the gluings between $E_{i-1,Q_{i-1}}$
and $E_{i,P_i}$ to a more general gluing would also decrease the
vanishing at $Q_i$. Therefore, this family is not part of a family
of limit linear series of dimension larger than $\rho$. It follows
that it can be deformed to the generic curve and will be part of a
family of dimension $\rho$ on the generic curve.

It remains to show that no coherent subsystem $(E',V')$
contradicts stability for any value of $\alpha$ or equivalently
that ${k'\over r'}\le {k\over r}$. Assume by contradiction that
this subsystem existed. Then the restriction to the curve $C_1$
would give a coherent subsystem of the restriction of $(E,V)$.
Note that $V$ is a generic subspace of sections of
$E(-(d_1-k_1-1)Q)$. The latter is a generic vector bundle of
degree $r(k_1+1)+d_2>k$. Then, Thm. 5.4 of \cite{LN} applies and
such a $V'$ cannot exist.

\bigskip
\bigskip
Consider now the case in which $0\not= d_2\ge k_2$ and assume that
$g+k_1-d_1\ge 2$. The latter condition is equivalent to the fact
that the generic vector bundle of rank $r$ and degree $d$ has
fewer than $k$ sections.

 On $C_1$ take the vector bundle to be the direct sum of $h$
 generic vector bundles of rank $\bar r$ and degree $\bar d$.

  On the curve $C_2$ take
the vector bundle
$$ {\mathcal O}(d_1Q_2)^{k_2}\oplus
L^2_1\oplus...\oplus L^2_{r-k_2}$$ where the $L^2_i$ are generic
line bundles of degree $d_1$. Take the gluing at $P_2$ so that the
subbundle  $ {\mathcal O}(d_1Q_2)^{k_2}$ glues inside the fiber at
$Q_1$ of  the space of dimension $d_2$ of sections of $E_1$ with
maximum vanishing $d_1$ at $Q_1$. The fiber of
$L^2_1\oplus...\oplus L^2_{r-k_2}$ at $P_2$ intersects
$H^0(E_1(-d_1Q_1))$ in a vector space of dimension $d_2-k_2$.
Consider the space of sections of $H^0(E_2(-(d_1-1)Q_2))$ that
glue with this space. They generate a space of dimension $d_2-k_2$
at $Q_2$. Denote this by $\bar V$. This is a space of dimension
$k_2$.

Define $W_2=H^0(E_2(-d_1Q_2))$.

 On the curve $C_i,\ i=3,...,k_1+2$ take
the vector bundle
$$ {\mathcal O}((2i-4)P+(d_1-2i+4)Q)^{r+k_2-d_2}\oplus
L^i_1\oplus...\oplus L^i_{d_2-k_2}$$ where the $L^i_j$ are generic
line bundles of degree $d_1$. At $P_3$ take the gluing so that
$L^3_1\oplus...\oplus L^3_{d_2-k_2}$ on the curve $C_3$ glues with
$\bar V$ defined above.

For $i\ge 3$, take the gluing so that the subbundle
$L^i_1\oplus...\oplus L^i_{d_2-k_2}$ on the curve $C_i$ glues with
the corresponding subbundle on the curve $C_{i+1}$.

Define $W_i$ by induction as the space of sections of
$H^0(E_i(-(i-3)P-(d_1-i+2)Q)$ that glues with $W_{i-1}$ at $P_i$.

On the curves
$$C_{(k_1+1)\alpha +\beta +1 },\ \alpha =1,...,g+k_1-d_1-2, \
\beta =1,...,k_1+1$$ take the following vector bundles:

 If $\beta =1$ take
$$ {\mathcal O}((k_1\alpha +1)P+(d_1-k_1\alpha -1)Q)^{k_2}\oplus
L^i_1\oplus...\oplus L^i_{r-k_2}$$

For $\beta =2,...k_1+1$ take
$$ {\mathcal O}((k_1\alpha +2\beta -2)P+(d_1-k_1\alpha -2\beta
+2)Q)^{r}.$$

For $\beta =2,..., k_1+1$, the gluing is generic. For $\beta =1,
\alpha =1$ take the gluing so that $ {\mathcal O}((k_1\alpha
)P+(d_1-k_1\alpha  )Q)^{k_2}$ glues with $W_{k_1+2}$.

Define $W_i$ inductively using the space of sections
$$H^0(E_i(-[k_1\alpha P+(d_1-k_1\alpha)Q])),\ \beta =1$$
$$H^0(E_i(-[(k_1\alpha +\beta -2)P+(d_1-k_1\alpha -2)Q])),\ \beta >1$$
that glues with the previous $W_{i-1}$.

For $\beta =1, \alpha>1$, glue the fiber at $P_i,\ i=(k_1+1)\alpha
+2$ of ${\mathcal O}((k_1\alpha +1)P+(d_1-k_1\alpha -1)Q)^{k_2}$
with the fiber at $Q_{i-1}$ of $W_{i-1}$.

On the remaining curves, the vector bundles are direct sums of
generic line bundles of degree $d_1$ and the gluings are generic.

 Define $t=g+k_1-d_1-1$. One can again define  $W_i$  inductively using the space
of sections of
$$H^0(E_i(-[(i-t-2) P+(d_1-i+t+1)Q])).$$
that glue with the previous $W_{i-1}$.

The resulting vector bundle is stable  because the restriction of
the vector bundle to each component is semistable and the
destabilizing subbundles in the various components do not glue
with each other.

In order to compute the dimension of the set of such vector
bundles, we need to add the dimensions of the families of the
vector bundles on each component, the dimensions of the families
of gluings, subtract the dimensions of the families of
automorphisms of the vector bundle restricted to each component
and finally add one as the resulting bundle being stable, it has a
one dimensional family of automorphisms.

The first bundle moves in an $h$ dimensional family and has an $h$
dimensional family of automorphisms.

The second one moves in an $[r-k_2]$-dimensional family and has a
$[(k_2)^2+r-k_2]$-dimensional family of automorphisms.

Each of the next $k_1$  bundles moves in a $[d_2-k_2]$-dimensional
family and has a $[(r+k_2-d_2)^2+d_2-k_2]$-dimensional family of
automorphisms. The first gluing varies in a
$[d_2k_2+r(r-k_2)]$-dimensional family. The next $k_1$ gluings
vary in  $[(d_2-k_2)^2+r(r-d_2+k_2)]$-dimensional families.

 On the curves corresponding to $\beta =1$, the vector bundle
 varies in a $[r-k_2]$-dimensional family and has a
$[(k_2)^2+r-k_2]$-dimensional family of automorphisms while the
gluing depends on $[k_2^2+r(r-k_2)]$ parameters.

On the curves corresponding to $\beta >1$, the vector bundles are
completely determined, they have an $r^2$-dimensional family of
endomorphisms and the gluings are free.

On the remaining curves, the vector bundles depend on
$r$-parameters, have an $r$-dimensional family of endomorphisms
and the gluings are free.

Therefore, the dimension of the family is
$$[h-h]+[r-k_2-(k_2)^2-(r-k_2)
+d_2k_2+r(r-k_2)]$$
$$+k_1[d_2-k_2-(r-(d_2-k_2))^2-(d_2-k_2)
+(d_2-k_2)^2+r(r-d_2+k_2)]$$
$$(g+k_1-d_1-2)[r-k_2-k_2^2-(r-k_2)+k_2^2+r(r-k_2)]+(g+k_1-d_1-2)k_1(-r^2+r^2)$$
$$+(g-(k_1+1)(g+k_1-d_1-1)-1)(r-r+r^2)+1=\rho$$

We need to define the spaces of sections on every component of the
curve and check that this gives rise  to only one limit series for
each bundle as defined above.

Define
$$W'_{k_1+2}=H^0(E_{k_1+2}(-2k_1P-(d_1-2k_1)Q)).$$
By definition, $$E_{k_1+2}=({\mathcal
O}((2k_1)P+(d_1-2k_1)Q))^{r+k_2-d_2}\oplus L^{k_1+2}_1\oplus
...\oplus L^{k_1+2}_{d_2-k_2}$$ where the last $d_2-k_2$ terms are
generic line bundles of degree $d_1$. Therefore,
$E_{k_1+2}(-(2k_1)P-(d_1-2k_1)Q)$ consists of a direct sum of
$r+k_2-d_2$ copies of the trivial line bundle and $d_2-k_2$
generic line bundles of degree zero. It follows that $W'_{k_1+2}$
has dimension $r+k_2-d_2$.

For $i<k_1+2$, define $W'_i$ by descending induction as the space
of sections of $E_i(-(i+k_1-2)P-(d_1-i-k_1+1)Q)$ that glues with
$W'_{i+1}$ at $Q_i$.

Define now the spaces of sections of the limit linear series.

For $i=1$, the space of sections is
$$H^0(E_1(-(d_1-k_1+1)P))\oplus W'_1.$$
 The vanishings of the sections of this space at $P_1$ are
$$0,1,k_1-1, k_1$$ where each vanishing is repeated $r$ times
except for the last one which is repeated $k_2$ times. The
vanishing of the sections of this space at $Q_1$ are
$$d_1,d_1-1...,d_1-k_1+1, d_1-k_1$$ where each vanishing is repeated $r$ times
except for the first one which is repeated $d_2$ times and the
last one which is repeated $r+k_2-d_2$ times.

For $i=2$, the space of sections is
$$W_2\oplus \bar V\oplus H^0(E_2(-P_2-(d_1-k_1)Q_2))\oplus
W'_2$$ The vanishing of this space of sections at $P_2$ is
$0,1,...k_1-1,k_1$, where each vanishing is repeated $r$ times
except for the first one that appears $d_2$ times and the last one
that appears $r-d_2+k_2$ times. At $Q_2$, the vanishing is $d_1,
d_1-1,..., d_1-k_1-1$ where each vanishing appears $r$ times
except for the first one which appears $k_2$ times, the second one
which appears $d_2-k_2$ times and the last one that appears
$r+k_2-d_2$ times.

For $i=3,...,k_1+2$, the space of sections is
$$W_i\oplus H^0(E_i(-(i-2)P-(d_1-2i+4)Q))\oplus
H^0(E_i(-(2i-3)P-(d_1-k_1-i+2)Q))$$ $$ \oplus
H^0(L^i_1\oplus...\oplus
L^i_{d_2-k_2}(-(2i-4)P-(d_1-2i+3)Q))\oplus W'_i.$$ The vanishings
of the sections of this space at $P_i$ are
$$i-3,i-2,..2i-6, 2i-5, 2i-4 ...k_1+i-3, k_1+i-2$$ where each vanishing is repeated $r$ times
except for the first one which is repeated $k_2$ times, the
vanishing $2i-5$ which is repeated $d_2-k_2$ times and the last
one which is repeated $r+k_2-d_2$ times. The vanishings of the
sections of this space at $Q_i$ are
$$d_1-i+2,d_1-i+1...d_1-2i+5, d_1-2i+4, d_1-2i+3...,d_1-k_1-i+2, d_1-k_1-i+1$$
where each vanishing is repeated $r$ times except for the first
one which is repeated $k_2$ times, the vanishing $d_1-2i+3$ which
is repeated $d_2-k_2$ times and the last one which is repeated
$r+k_2-d_2$ times.

For $i=\alpha (k_1+1)+\beta$ the space of sections is
$$W_i \oplus H^0(E_i(-(k_1\alpha +\beta -1)P-(d_1-(k_1\alpha +2\beta -2)Q))\oplus $$
$$ \oplus H^0(E_i(-(k_1\alpha +2\beta -1)P-(d_1-(k_1(\alpha +1)+\beta )Q)).$$
The vanishings of the sections of this space at $P_i$ are
$$i-\alpha-3,..k_1\alpha+2\beta -4,
k_1\alpha+2\beta -2,k_1\alpha+2\beta -1,...k_1(\alpha+1)+\beta
-1$$ where each vanishing is repeated $r$ times except for the
first one which is repeated $k_2$ times. The vanishings of the
sections of this space at $Q_i$ are
$$d_1-i+\alpha+2,..d_1-k_1\alpha-2\beta +3,
d_1-k_1\alpha-2\beta +2,d_1-k_1\alpha-2\beta
,...d_1-k_1(\alpha+1)-\beta $$ where each vanishing is repeated
$r$ times except for the first one which is repeated $k_2$ times.

Write $t=g+k_1-d_1-1$. then for $i>t(k_1+1)+1$, the space of
sections is $$W_i\oplus H^0(E_i(-(i-t-1)P-(d_1-i+t-k_1+1)Q))$$ The
vanishings of the sections of this space at $P_i$ are
$$i-t-2,i-t-1,...i-t+k_1-2$$
where each vanishing is repeated $r$ times except for the first
one which is repeated $k_2$ times. The vanishings of the sections
of this space at $Q_i$ are
$$d_1-i+t+1,d_1-i+t,..d_1-i+t-k_1+1$$
where each vanishing is repeated $r$ times except for the first
one which is repeated $k_2$ times.

As in the case $d_2<k_2$, this defines a limit linear series with
$b=d_1$ which is not part of a larger family of limit linear
series, therefore it deforms to the generic curve.

\bigskip
It only remains to prove the claim that for every coherent
subsystem $(E',V')$ of rank $r'$, degree $d'$ and with  $dim
V'=k'$, $k'/r'\le k/r$. Suppose by contradiction that there exists
a subsystem with $k'/r'> k/r$. Note that by our choice of gluing,
the space of sections that we take on the first vector bundle
$E_1$ is in fact a generic subspace of dimension $k$ of
$H^0(E_1(-(d_1-k_1)Q_1)$. Then if $k_2<d_2$ or $k_2=d_2, h=1$ from
\cite{LN} Th 5.4, the restriction of $(E',V')$ to $C_1$ cannot
exist.

It remains to deal with the case $k_2=d_2, h\not= 1$. Then
$E_1=F_1\oplus ...\oplus F_h$ where each $F_j$ is an
indecomposable bundle of rank $\bar r$ and degree $\bar d=\bar
rd_1+\bar d_2, \ k_2=d_2=h\bar d_2$.
 As $E_1$ is semistable,  the slope of $E'$ is at most the
slope of $E_1$. Again, $V'$ is a subspace of the space of sections
of $E'(-(d_1-k_1)Q_1)$ . The latter is a vector bundle of rank
$r'$ and degree $d'-r'(d_1-k_1)=r'k_1+d'_2$.  Then, $k'\le
r'k_1+d'_2$. Hence,
$$ {k'\over r'}\le
{r'k_1+d'_2\over r'}\le {rk_1+d_2\over r}={rk_1+k_2\over
r}={k\over r}.$$

\bigskip

When $d_2=k_2=0$, take on $C_1$ the direct sum of $r$ generic line
bundles of degree $d_1$.  On
$$C_{\alpha k_1+\beta +1},\ \alpha=0,...,g-2-d_1+k_1,\ \beta =1,...,k_1$$ take
$${\mathcal O}((2\beta +\alpha (k_1-1)-1)P+(d_1-2\beta -\alpha
(k_1-1)+1)Q)^r$$ On the remaining curves take direct sums of
generic vector bundles of degree $d_1$.

As spaces of sections of the limit linear series take on
$C_{\alpha k_1+\beta +1}$
$$H^0(E_i(-(i-2-\alpha)P-(d_1-2\beta-\alpha (k_1-1)+1)Q))\oplus$$
$$\oplus H^0(E_i(-(2\beta+\alpha (k_1-1))P-(d_1-i-k_1+\alpha +1)Q)$$

Write $t=g-d_1+k_1-1$. For $i\ge k_1t+2$, take as space of
sections $H^0(E_i(-(i-t-1)P-(d_1-i+t-k_1+1))Q))$. We leave it to
the reader to check that one gets a limit linear series with the
value of the parameter $b$ in the definition being $d_1$ and the
dimension of the family being $\rho$.

\bigskip

Consider now the case in which $d+r(1-g)\ge k$. Then, the vector
bundle is going to be generic while the limit linear series for a
fixed bundle will move in a space of dimension equal to the
dimension of the grassmannian $Gr(k, d+r(1-g))$. The details are
as follow:

On the curve $C_1$ take as before  a generic sum of $h$ vector
bundles of rank $\bar r$ and degree $\bar d$.
 On the curves
$C_2,...C_g$ take the vector bundle to be a direct sum of $r$
generic line bundles of degree $d_1$.

The integer $b$ that appears in III of the definition of limit
linear series above will be taken again to be $d_1$. Glue the
vector bundles $E_i$ to $E_{i+1}$ at the point $Q_i,P_{i+1}$ by
using generic gluing.

The vector bundle obtained in this way is stable because all the
restrictions to the elliptic components are semistable and  the
destabilising line bundles for each component do not glue with
each other (see \cite{T3,T4}). Moreover it corresponds to a
generic point of a component of the moduli space of vector bundles
of this rank and degree. Hence, it moves in a family of dimension
$r^2(g-1)+1$.

Take a space of sections of $E_1$ of dimension $k$ generic inside
 $E_1(-(g-1)Q_1)$. For each of the vanishings $g-1, g,
 ...,g+k_1-2$, this space contains $r$ independent sections with these
 orders of vanishing. It contains a $k_2$-dimensional space of
 sections with vanishing $g+k_1-1$.

 On $E_i$ one must then take $r$-dimensional spaces of sections
 with vanishing $g-i+j$ at $Q_i$ and $d_1-g+i-1-j$ at $P_i$,
 $j=0...k_1-1$ as well as a $k_2$ dimensional space of sections
 with vanishing $d_1-g+i-k_1-1$ at $P_i$ and $g-i+k_1$ at $Q_i$ that
 glues at $P_i$ with the corresponding space at $Q_{i-1}$.

 One can check that this gives rise to a limit linear series.
 The sum of the vanishings at $Q_{i-1}$ and $P_i$ of corresponding sections adds up
  to the minimum required ($d_1$) while the vanishing at $Q_g$ is
  $0, 1...k_1-1, k_1$ where each number
  is repeated $r$ times except for the last one which is repeated
  $d_2$ times. These are the minimum possible vanishings for the
  sections of a $k$ dimensional linear series at a point. Hence,
  the series is completely determined by the choice of the
  subspace of sections on the first curve. This first choice is
  equivalent to the choice of a $k$-dimensional vector space of
  $H^0(E_1(-(g-1)Q_1))$. The latter is a vector space of dimension
  $d-r(g-1)$. It follows that, for a fixed vector bundle, there is a
  $[k(d-r(g-1)-k)]$-dimensional family of limit linear series.
  Then, the set of coherent
  systems moves in a family of dimension $r^2(g-1)+1+k(d-r(g-1)-k)=\rho$.

Assume that the coherent system contained a coherent subsystem
corresponding to a subbundle $E'$ of rank $r'$ with $k'$ sections
where ${k'\over r'}>{k\over r}$. If $k<d+r(1-g)$ or $(r,d)=1$,
Theor.5.4 of \cite{LN} shows that the restriction to $C_1$ cannot
exist. When $k=d+r(1-g)$ and $(r,d)\not= 1$, the result follows
with arguments similar to the previous cases.

 This concludes the proof of the Theorem.

\end{section}


\bigskip




\begin{thebibliography}{cccc}
\bibitem [B]{B} L.Brambila-Paz {\it Non-emptiness of moduli spaces
of coherent systems} AG/0412285

\bibitem [BG]{BG} S.Bradlow, O.Garcia-Prada {\it An application of
coherent systems to a Brill-Noether problem} J.reine angew.Math.
{\bf 551} (2002), 123-143.

\bibitem [BGMN]{BGMN} S.Bradlow, O.Garcia-Prada, V.Munoz, P.Newstead
{\it Coherent systems and Brill-Noether theory}, Internat.J.Math.
{\bf 14}(2003), 683-733.

\bibitem [BGMMN] {BGMMN} S.Bradlow, O.Garcia-Prada,V.Mercat, V.Munoz, P.Newstead
{\it On the geometry of moduli spaces of coherent systems on
algebraic curves}, AG/0407523.

\bibitem [KN]{KN} A.King, P.Newstead, {\it Moduli of Brill-Noether
pairs on algebraic curves}, Internat J.Math {\bf 6}(1995),
733-748.

\bibitem [LP]{LP} J. Le Potier, {\it Faisceaux semistables et
systemes coherents} in "Vector bundles in algebraic geometry
(Durham 1993)", 179-239, LondonMath,Soc. lecture note series 1995.


\bibitem [LN]{LN} H.Lange, P.Newstead {\it Coherent systems on
elliptic curves} Intern. J.Math. {\bf 16} (2005), 787-805.

\bibitem [S]{S} N.Sundaram {\it Special divisors and vector
bundles} Tohoku Math J.{\bf 39}(1987), 175-213.


\bibitem [T1]{T1} M.Teixidor {\it Existence of vector bundles of
rank two with sections}, Adv.Geom. {\bf 5} (2005), 37-47.

\bibitem [T2]{T2} M.Teixidor {\it Brill Noether Theory for
stable vector bundles}, Duke Math. J. {\bf 62 \ N2}, (1991)
385-400

\bibitem [T3]{T3} M.Teixidor {\it Moduli spaces of semistable
vector bundles on tree-like curves}, Math.Ann. {\bf 290} (1991).
341-348.

\bibitem [T4]{T4} M.Teixidor {\it Moduli spaces of vector bundles on
reducible curves}, Amer J of Math. {\bf 219}(1995), 477-484.


\end{thebibliography}
\end{document}